\documentclass[12pt]{article} 
\input epsf 
\usepackage{amsmath} 
\def\Rea{\operatorname{Re}} 
\def\Ima{\operatorname{Im}} 
\def\grad{\mathrm{grad}} 
\def\C{{\mathbf{C}}} 
\def\R{{\mathbf{R}}} 
\def\s{\overline{s}} 
\begin{document} 
\title{On the number of solutions of a transcendental 
equation arising in the theory of gravitational lensing 
} 
\author{Walter Bergweiler\thanks{ 
Supported by the EU Research Training Network CODY, the  
ESF Networking Programme HCAA and the Deutsche Forschungsgemeinschaft, Be 1508/7-1.} 
\hskip.1in and Alexandre Eremenko\thanks{Supported by 
the NSF grant DMS-0555279 and by the Humboldt Foundation.}} 
\maketitle 
\begin{abstract} 
The equation in the title describes the number of bright images 
of a point source under 
lensing by an elliptic object with  
isothermal density. We prove that this equation 
has at most $6$ solutions. 
Any number of solutions 
from $1$ to $6$ can actually occur. 
\end{abstract} 
 
\begin{center} 
{\bf 1. Introduction} 
\end{center} 
\vspace{.1in}
 
We study the number of solutions of the equation 
\begin{equation}\label{1} 
z-\frac{k}{\sin\overline{z}}=w, 
\end{equation} 
in the region $|\Rea z|<\pi/2$ in the complex plane, where 
$w\in\C$ and $k>0$ are para\-meters. 
This equation is equivalent to the equation
\begin{equation}\label{1a} 
z=\arcsin\frac{k}{\overline{z}+\overline{w}},
\end{equation} 
which occurs as a model of gravitational 
lensing of a point source~$w$ by an elliptic 
object whose density equals $c/r$ on the homothetic ellipses $rE$ 
where $E$ is a fixed ellipse, cf.~\cite{H,Keeton,HL}. 
The branch of $\arcsin$ in (\ref{1a}) is the principal
branch which is defined in $\C\backslash[-1,1]$.
As in~\cite{H,HL} we assume that the density is
zero outside 
of $E$ and thus is equal
to $c/r$ on $rE$ only for $0<r\leq 1$. 
We note, however, that 
in the astronomy literature (cf., e.~g.,~\cite{Keeton})
it is usually assumed  
that this formula for the density holds for $0<r<\infty$;
see~\cite{HL} for a discussion of the different models.
We also refer to the above papers for the derivation of 
equations (\ref{1}) and~(\ref{1a}). 
 
The solutions of~(\ref{1}) that lie 
outside $E$ 
correspond to the so-called 
bright images of the source. 
Khavinson and Lundberg~\cite{HL} proved that  
the number of solutions of~(\ref{1}) in $|\Rea z|<\pi/2$ 
is finite and does not 
exceed~$8$.  
Up to $4$ bright images by a single lensing galaxy
have been observed by  
astronomers~\cite{HL,site}. 
 
In this paper we prove the following result. 
\vspace{.1in} 
 
\noindent{\bf Theorem.} {\em The number $n$ 
of solutions of equation  
$(\ref{1})$ in the region $|\Rea z|<\pi/2$ satisfies 
$n\leq 6$.} 
\vspace{.1in} 
 
Any number between $1$ and $6$ 
can actually occur. 
In fact, for each value of the para\-meter $k$ we will 
describe a partition of the $w$-plane into regions 
where the number of solutions is constant. 
This also yields that  
equation~(\ref{1}) has 
at least one solution for every $k>0$ and $w\in\C$,
but we do not include a formal proof of this fact.
 
The upper estimate in \cite{HL} is based on 
a miraculous trick: an application 
of Fatou's theorem from holomorphic dynamics. 
This method originated in~\cite{HS}, and it was applied 
by Khavinson and Neumann~\cite{HN} to obtain the estimate $5d-5$ 
for the number of solutions of the equation 
$$z-R(\overline{z})=0,$$ 
where $R$ is a rational function of degree~$d$. 
This equation 
with $$R(z)=\sum_{j=1}^d 
\frac{a_j}{z-z_j},\quad a_j>0,$$ 
describes gravitational lensing 
by $d$ coplanar point masses.  
The estimate $5d-5$ for the number of solutions is 
exact, even for this special form of~$R$; this 
was established by 
Rhie~\cite{Rhie,Rhie2}, 
who also conjectured the correct estimate $5d-5$. 
We refer to~\cite{HN1} for an exposition of the above 
and related results.
 
Our method does not use holomorphic dynamics. Instead it is based on 
elementary considerations from the 
theory of harmonic maps~\cite{Duren,Lyzzaik}. 
While our method is thus completely
different from that employed in~\cite{HL},
we shall use two facts established in~\cite{HL}.
One is that for all $w\in\C$ there are only finitely many 
solutions of~(\ref{1}) satisfying $|\Rea z|<\pi/2$, the other
one concerns the index of the image of the lines $|\Rea z|=\pi/2$
under a certain function~$f$; see section~2 below.
However, the proofs of these facts are elementary and do not
rely on dynamics.

For a specific elementary equation like~(\ref{1}), 
involving only two para\-meters $k$ and~$w$, the method 
we propose has the following advantage: 
it provides a method of determining the number 
of solutions for each value of the para\-meters. 
 
We hope to apply the same method 
to the 
more general equation 
\begin{equation}\label{shear} 
z=\arcsin\frac{k}{\overline{z}+\overline{w}}+\alpha{\overline{z}}, 
\quad \alpha\in\C, 
\end{equation} 
which describes gravitational lensing by an elliptic object 
of isothermal density and a shear 
\cite{Keeton,HL}. It is not clear how to use the arguments 
based on Fatou's theorem in the case that $\alpha\neq 0$. 
We will further discuss this at the end of 
the paper. 
 
We thank A.\ Gabrielov, C.\ R.\ Keeton, D.\ Khavinson,
E.\ Lundberg, Sun Hong Rhie and two referees for many useful
remarks and suggestions.
\vspace{.2in} 
 
\begin{center} 
{\bf 2. The harmonic map $f$ and its caustic} 
\end{center} 
\vspace{.1in}
 
We are solving the equation $f(z)=w$ where 
\begin{equation} \label{deff}
f(z)=z+\overline{g(z)},\quad g(z)=-\frac{k}{\sin z}, 
\end{equation} 
and $z\in D^0=\{ z\in\C:|\Rea z|<\pi/2\}$. 
The Jacobian determinant of $f$ is 
$$J(z)=1-|g'(\overline{z})|^2=1-|g'(z)|^2.$$ 
Our map $f$ is smooth and finite, 
which means that every point has only finitely many
preimages in $D^0$.
The last property (established in~\cite{HL}) follows from
the fact that solutions of the equation
$f(z)=w$ are fixed points of the function
$g_w=w-\overline{g}$,
so they also satisfy the equation
$z=(g_w\circ g_w)(z)$ where the right hand side is analytic.
Thus the $w$-points of $f$ are isolated, and since $f(z)\to\infty$
 as $z\to\infty$ in $D^0$, there are only finitely many  $w$-points 
in $D^0$.

Harmonic maps with discrete preimages of points
are called ``light'' in~\cite{Lyzzaik}.  
Our map $f$ preserves the orientation in the open set 
$D^+=\{z\in \C: J(z)>0\}$ and reverses the 
orientation in the complementary open set $D^-=\{z\in \C: J(z)<0\}$. 
The common boundary of $D^+$ and $D^-$ with respect to $D^0$ 
is given by the equation 
$$|g'(z)|= 
k\left|\frac{\cos z}{\sin^2 z}\right|=1.$$ 
We call this common boundary $\gamma$. 
In the astronomy literature, $\gamma$ is called the 
{\em critical curve}, 
and its image $\Gamma=f(\gamma)$ is called the {\em caustic}. 
 
\begin{figure} 
\epsfxsize=3.2in 
\centerline{\epsffile{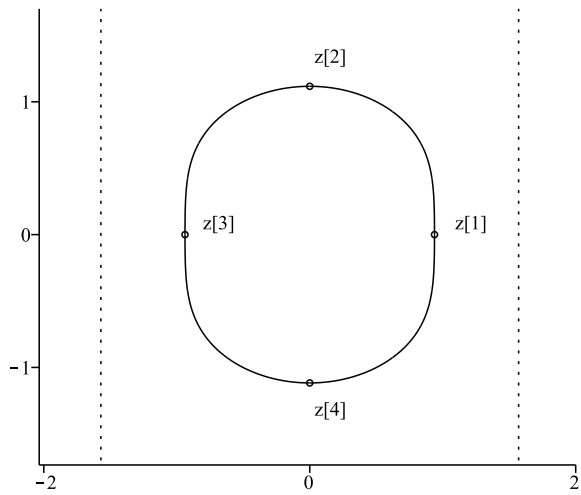}} 
\caption{Critical curve for $k=1.1$.} 
\label{f1} 
\end{figure}

\begin{figure} 
\epsfxsize=3.2in 
\centerline{\epsffile{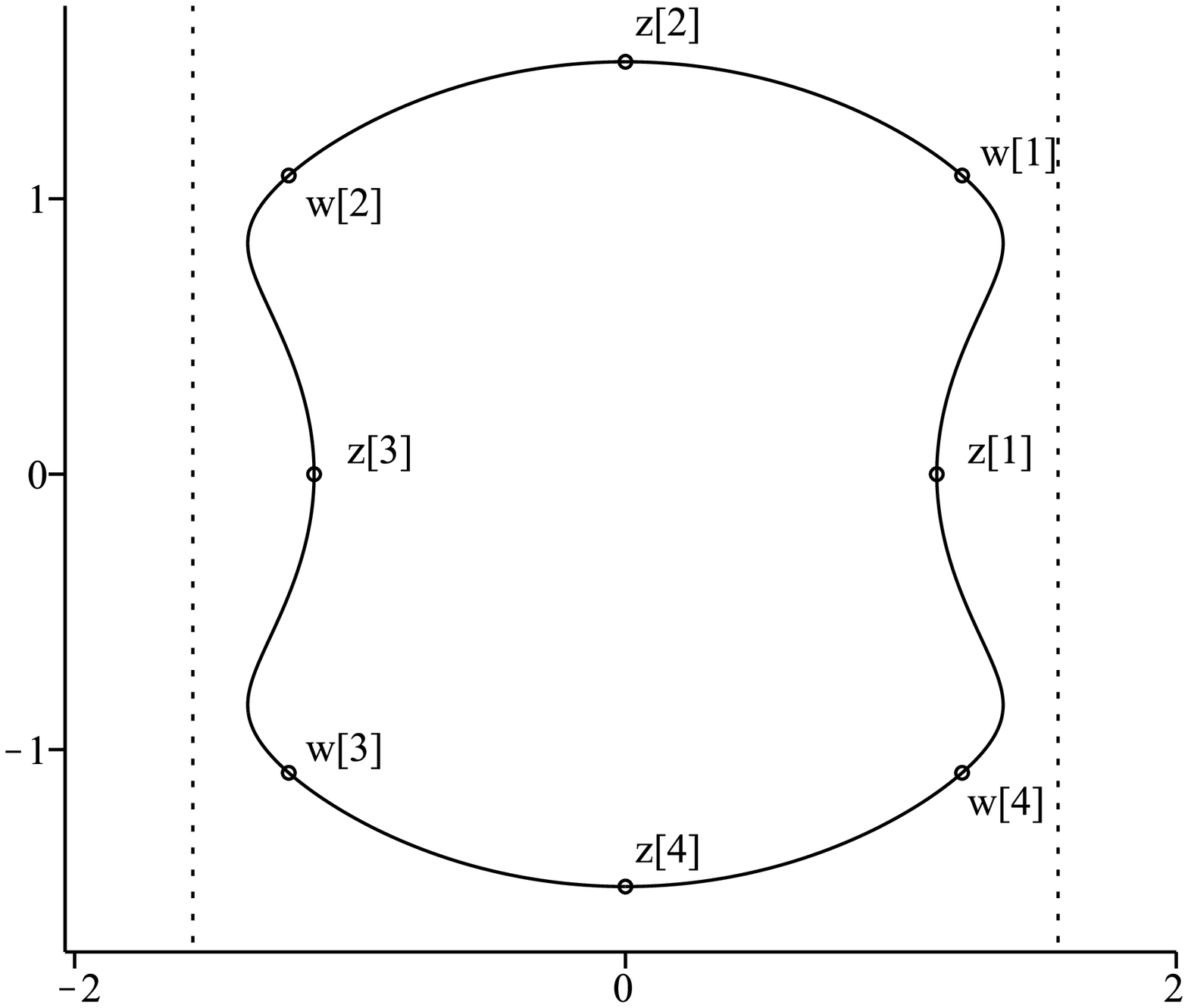}} 
\caption{Critical curve for $k=1.92$.} 
\label{f2} 
\end{figure}

\begin{figure} 
\epsfxsize=3.2in 
\centerline{\epsffile{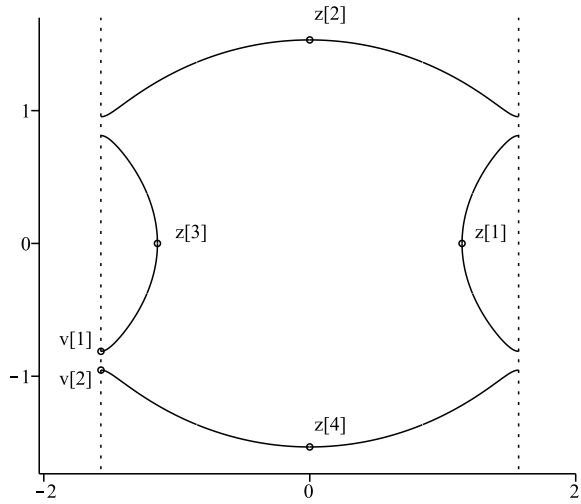}} 
\caption{Critical curve for $k=2.01$.} 
\label{f3} 
\end{figure}

For small positive~$k$, the critical curve $\gamma$ consists of 
a single smooth Jordan curve surrounding the pole at~$0$. 
This picture persists for $0<k<2$. At $k=2$, the critical curve $\gamma$ 
bifurcates into four smooth simple curves with endpoints 
on $\partial D^0$; see Figures \ref{f1}--\ref{f3}. 
 
Our count of the number of solutions is based 
on the Argument Principle for light mappings;  
see~\cite{Cristea,DHL,Lyzzaik}.  
Let $D$ be a region 
bounded by finitely many disjoint 
Jordan curves on the Riemann sphere. 
We para\-metrize the boundary curves so 
that the region stays on the 
left. Assuming that a smooth map $f:D\to\overline{\C}$, 
continuous in the closure of~$D$, 
never takes the values $w$ and $\infty$ 
 on the boundary $\partial D$, 
and that $J(z)\neq 0$ in~$D$, 
let $N$ and $P$ be the numbers of $w$-points and poles 
of~$f$, respectively. Then 
\begin{equation}\label{fo} 
N-P=\pm I_{w}(f(\partial D)), 
\end{equation} 
where $I_{w}(f(\partial D))$ is the index (or winding number) 
of the curves $f(\partial D)$ about~$w$. Thus 
$$I_w(f(\partial D))=\frac{1}{2\pi i}\int_{f(\partial D)}\frac{d\zeta}{\zeta-w}= 
\frac{1}{2\pi i}\int_{\partial D}\frac{ 
df(z)}{f(z)-w},$$ 
where 
$$df=f_zdz+f_{\overline{z}}d\overline{z}.$$ 
If $J(z)>0$ in $D$ so that 
the map preserves the orientation, we choose the plus sign 
in~(\ref{fo}), and if it reverses the orientation, we 
choose the minus sign. 
 
Our Theorem is an immediate consequence of the following 
propositions. 
 
\vspace{.1in} 
\noindent 
{\bf Proposition 1.} {\em If $k>0$, 
then $|I_w(f(\partial D^-))|\leq 2$ for all $w\in\C$. } 
\vspace{.1in} 

Here and in what follows the expression ``for all $w$'' 
means ``for all $w$ for which the index is defined''. 
A more careful analysis, which we do not 
include in this paper, would 
show that $I_w(f(\partial D^-)) 
\leq 0$ for all $k>0$ and all~$w$.  
This would imply that equation~(\ref{1}) has at 
least one orientation-reversing solution for all 
$k>0$ and $w\in\C$. 
\vspace{.1in} 
 
\noindent 
{\bf Proposition 2.} {\em 
If $k\geq 2$, then $|I_w(f(\partial D^+))|\leq 3$ for all 
$w\in\C$.} 
\vspace{.1in} 
 
In these propositions, the boundaries 
$\partial D^-$ and $\partial D^+$ 
are understood with respect to 
the extended plane; 
if $0<k<2$ then $\partial D^-=\gamma$, while if 
$k>2$, then $\partial D^-$ consists of four components of 
$\gamma$ and four vertical intervals. The open set $D^+$ is  
always unbounded. To show that the Argument Principle
applies to $D^+$ we exhaust $D^+$
 by regions of the form $D^+\cap\{ z:|\Ima z|<R\}$
with $R\to\infty$.
  
Proposition 2 is actually true for every $k>0$. This follows from 
the argument in \cite{HL} using Fatou's theorem. We give a proof independent 
of Fatou's theorem 
for the case that $k\geq 2$, which suffices for our purposes. 
 
To derive our Theorem from the propositions, 
we first assume that $w\not\in\Gamma$ and
consider two cases. 
 
If $0<k<2$, we apply the Argument Principle and Proposition 1 
to $D^-$, which contains one pole, and obtain 
that the number of orientation-reversing 
solutions of~(\ref{1}) is at most~$3$. 
Then we apply the Argument Principle to~$D^+$. 
The boundary of $D^+$ consists of the curve $\gamma$ and 
two vertical lines. The image of the two vertical lines 
is easy to study and its index about 
any point in the plane has absolute 
value at most~$1$, a fact established in~\cite{HL}. 
Thus $|I_w(f(\partial D^+))|\leq 3$ by Proposition~1. 
Since there are no poles in $D^+$, 
this implies that 
the number of orientation-preserving solutions of~(\ref{1})  
is at most~$3$. 
Thus there are at most $6$ solutions in this case. 
 
If $k\geq 2$, the argument is similar. By Proposition 2, 
the number 
of orientation-preserving solutions is at most~$3$, 
and by Proposition 1, the number of orientation-reversing solutions 
is at most $2+1=3$. 
So our equation 
has at most 
$6$ solutions in this case as well. 
 
If $w\in\Gamma$ we apply the following general fact.
\vspace{.1in}

\noindent
{\bf Proposition 3.} {\em Let $f:D\to\C$ be a harmonic map
defined in a region  $D$ in $\C$. Suppose that every $w\in\C$ 
has at most $m$ preimages, where $m<\infty$.
Then the set of points which have $m$ preimages is open.}
\vspace{.1in}

This is not true for arbitrary smooth maps.
As we found no reference for Proposition 3, we include a
proof in section 5.

This completes the 
derivation of our Theorem from Propositions 1--3. 
\vspace{.1in} 
 
\begin{figure} 
\epsfxsize=3.2in 
\centerline{\epsffile{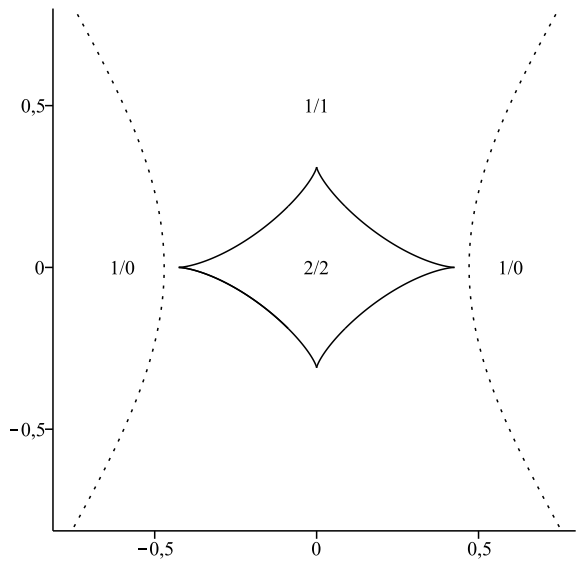}} 
\caption{The caustic $\Gamma$  and, in dotted lines,
 $f(\partial D^0)$ for $k=1.1$. Here
$m/n$ indicates the number of
orientation-reversing/preserving solutions.}
\label{f4} 
\end{figure}

\begin{figure} 
\epsfxsize=3.2in 
\centerline{\epsffile{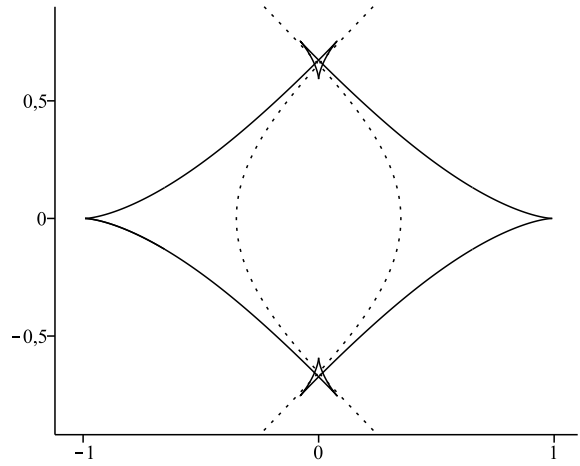}} 
\caption{The caustic and the image of the boundary of $D^0$ for $k=1.92$.} 
\label{f5} 
\end{figure}

\begin{figure} 
\epsfxsize=3.2in 
\centerline{\epsffile{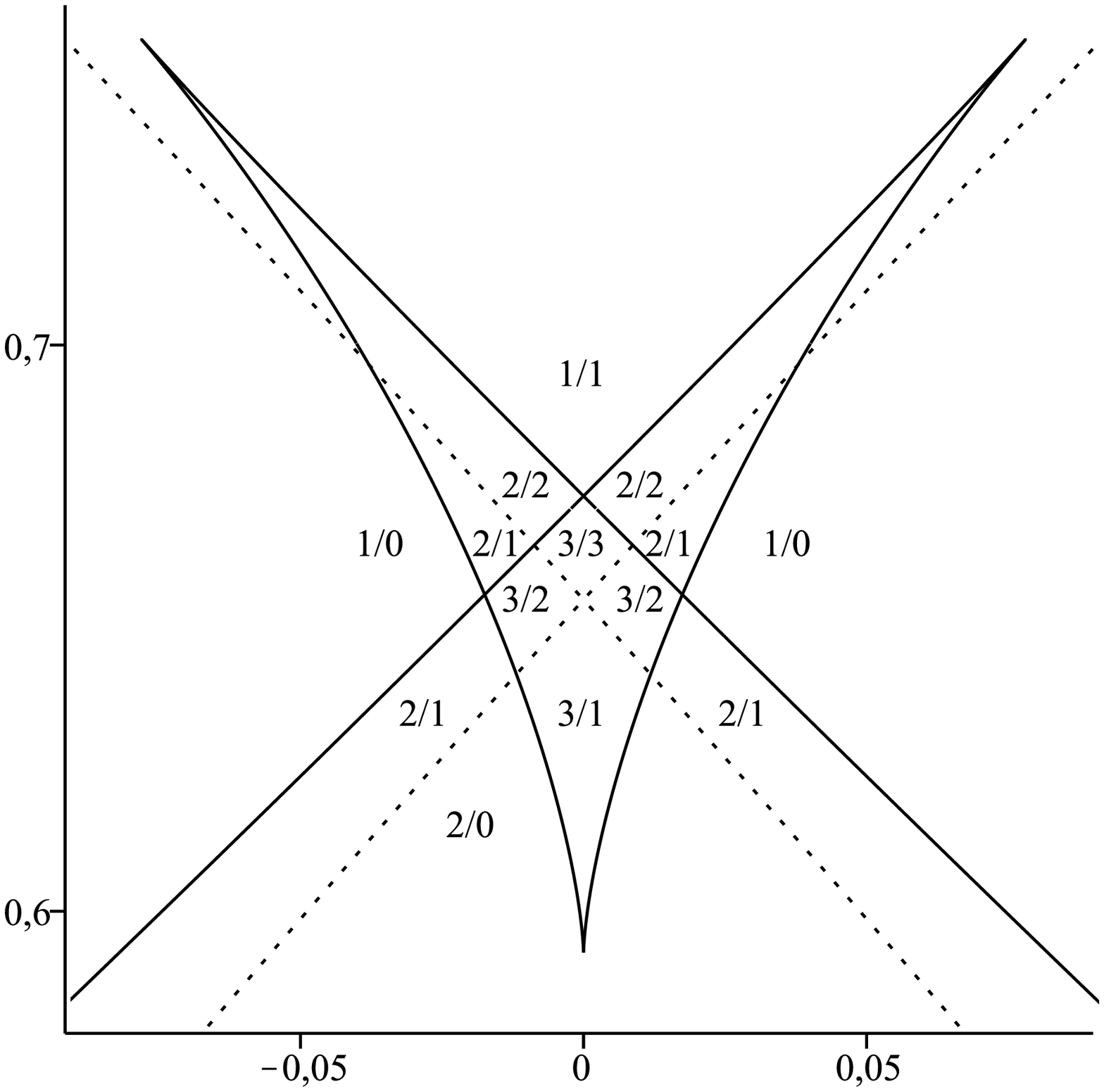}} 
\caption{Magnification of detail from Figure \ref{f5}.
Here $m/n$ indicates the number of
orientation-reversing/preserving solutions.}
\label{f6} 
\end{figure}

\begin{figure} 
\epsfxsize=3.2in 
\centerline{\epsffile{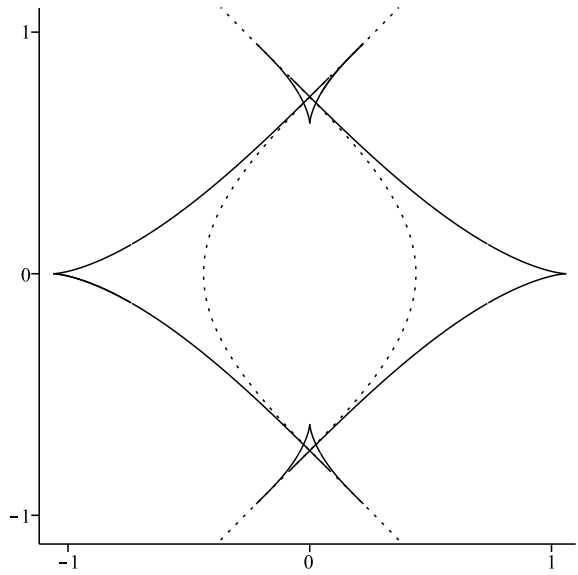}} 
\caption{The caustic and the image of the boundary of $D^0$ for $k=2.01$.} 
\label{f7} 
\end{figure}

\begin{figure} 
\epsfxsize=3.2in 
\centerline{\epsffile{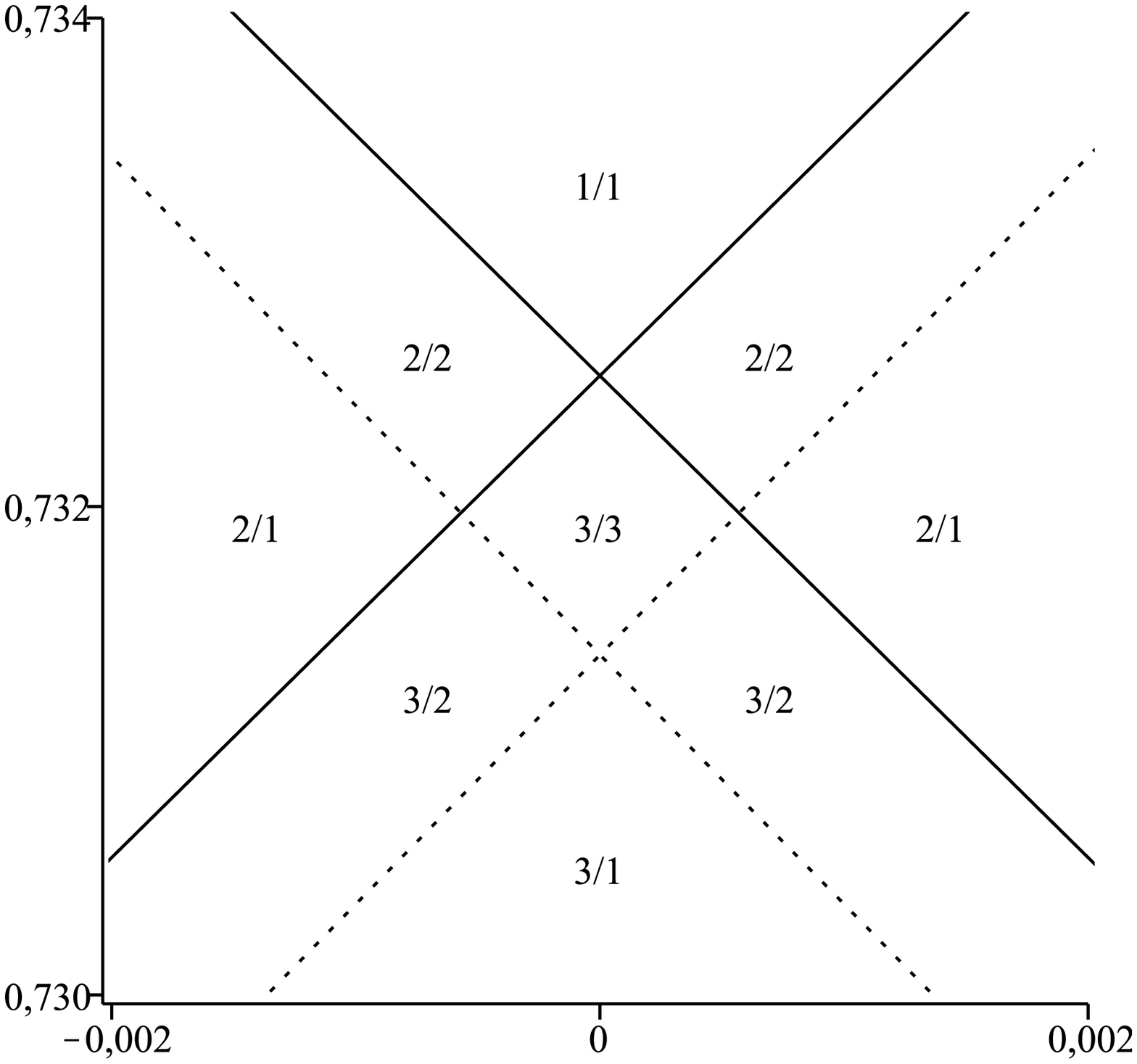}} 
\caption{Magnification of detail from Figure \ref{f7}.
Here $m/n$ indicates the number of
orientation-reversing/preserving solutions.}
\label{f7a} 
\end{figure}

Figures~\ref{f4}--\ref{f7a} show  
the images of the boundaries $f(\partial D^+)$ 
and $f(\partial D^-)$. 
The numbers of solutions of~(\ref{1}) are written in 
the regions complementary to these images. 
The notation $m/n$ in Figures~\ref{f4}, \ref{f6} and~\ref{f7a} 
means that for $w$ in the indicated region there are $m$ 
orientation-reversing and $n$ orientation-preserving 
solutions. 
 
The essential bifurcation occurs at the point 
$k=2/\sqrt{3}\approx 1.1546$. 
 In particular, $5$ or $6$ solutions 
are only possible for 
$$\frac{2}{\sqrt{3}}<k<k_0
=\frac{\pi^2}{2\sqrt{\pi^2-4}}\approx 2.0368$$
and the region 
in the $w$-plane where the number of solutions is $5$ or $6$ 
is rather small. Perhaps this explains the fact that 
$5$ or $6$ bright images in a single elliptic lens
do not seem to have been observed 
by astronomers yet. 
 
For some readers these pictures produced by Maple 
will be sufficiently convincing; these readers may skip 
the next three sections 
and pass to the remarks in the end of the paper. 
The formal proofs which we give below 
show that the pictures actually represent correctly all 
essential features of our curves, necessary to determine 
their indices about every point in the plane. 
As we mentioned before, some of those features are rather 
small, and one has to be sure that nothing was missed on 
a still smaller scale. 
\vspace{.2in} 

\begin{center}  
{\bf 3. The cusps of the caustic}  
\nopagebreak[4]  
\end{center}  
\nopagebreak[4]  
\vspace{.1in}  
\nopagebreak[4] 
The curve $\gamma$ is given by the equation 
\begin{equation} 
\label{11} 
|g'(z)|=k\left|\frac{\cos z}{\sin^2 z}\right|=1. 
\end{equation} 
We use the local theory of harmonic mappings following the 
paper by Lyzzaik~\cite{Lyzzaik}. 
We mention that Lyzzaik considered only harmonic maps in
simply connected domains while our map $f$ has a pole.
However, since the results are local, this does not affect 
the arguments.
We para\-metrize the curve by  
\begin{equation} 
\label{12} 
t=-\arg g'. 
\end{equation} 
This corresponds to the {\em counterclockwise} 
motion on the curve $\gamma$. 
First we determine the cusps of $f(\gamma)$. 
Let $z(t)$ be the para\-metrization of $\gamma$. 
The cusps in the image 
will occur when  
\begin{equation}\label{dfz} 
\frac{d}{dt}f(z(t))=0. 
\end{equation} 
A simple computation in \cite[(2.4), (2.5)]{Lyzzaik} 
shows that~(\ref{dfz}) yields 
\begin{equation}\label{L} 
\Rea(z'(t)e^{-it/2})=0. 
\end{equation} 
In order to compute the argument of $z'(t)$ we note that 
the curve $\gamma$ which is para\-metrized by $z(t)$ is a level curve  
of $\log |g'|$. This implies that 
$$\arg z'(t)=\arg\grad(\log|g'(z(t))|)-\frac{\pi}{2}=\arg\left( 
\frac{g'(z(t))}{g^{\prime\prime}(z(t))}\right)-\frac{\pi}{2}.$$ 
So, using~(\ref{12}), 
$$\arg \left(z'(t)e^{-it/2}\right)=-\arg g^{\prime\prime}(z(t))+ 
\frac32\arg g'(z(t))-\frac{\pi}{2},$$ and 
condition~(\ref{L}) becomes 
\begin{equation}\label{13} 
\frac{ 
g^{\prime\prime}(z)^2}{g'(z)^3}=\frac{(1+\cos^2 z)^2}{k \cos^3 z}>0. 
\end{equation} 
So to locate the cusps we need to solve two simultaneous 
equations~(\ref{11}) and~(\ref{13}), that is, 
$$k\left|\frac{\cos z}{\sin^2 z}\right|=1,\quad 
\frac{(1+\cos^2 z)^2}{\cos^3 z}>0,$$ 
where we used that $k>0$ to drop $k$ from 
the second equation. 
Putting $s=\cos z$ we obtain the algebraic equations 
\begin{equation} 
\label{I} 
k^2\frac{s\s}{(1-s^2)(1-\s^2)}=1 
\end{equation} 
and 
\begin{equation} 
\label{II} 
\frac{(1+s^2)^2}{s^3}-\frac{(1+\s^2)^2}{\s^3}=0. 
\end{equation} 
The second equation expresses the condition 
$(g^{\prime\prime})^2/(g')^3\in \R$; we will later select 
those solutions that satisfy~(\ref{13}). 

\begin{figure} 
\epsfxsize=3.2in 
\centerline{\epsffile{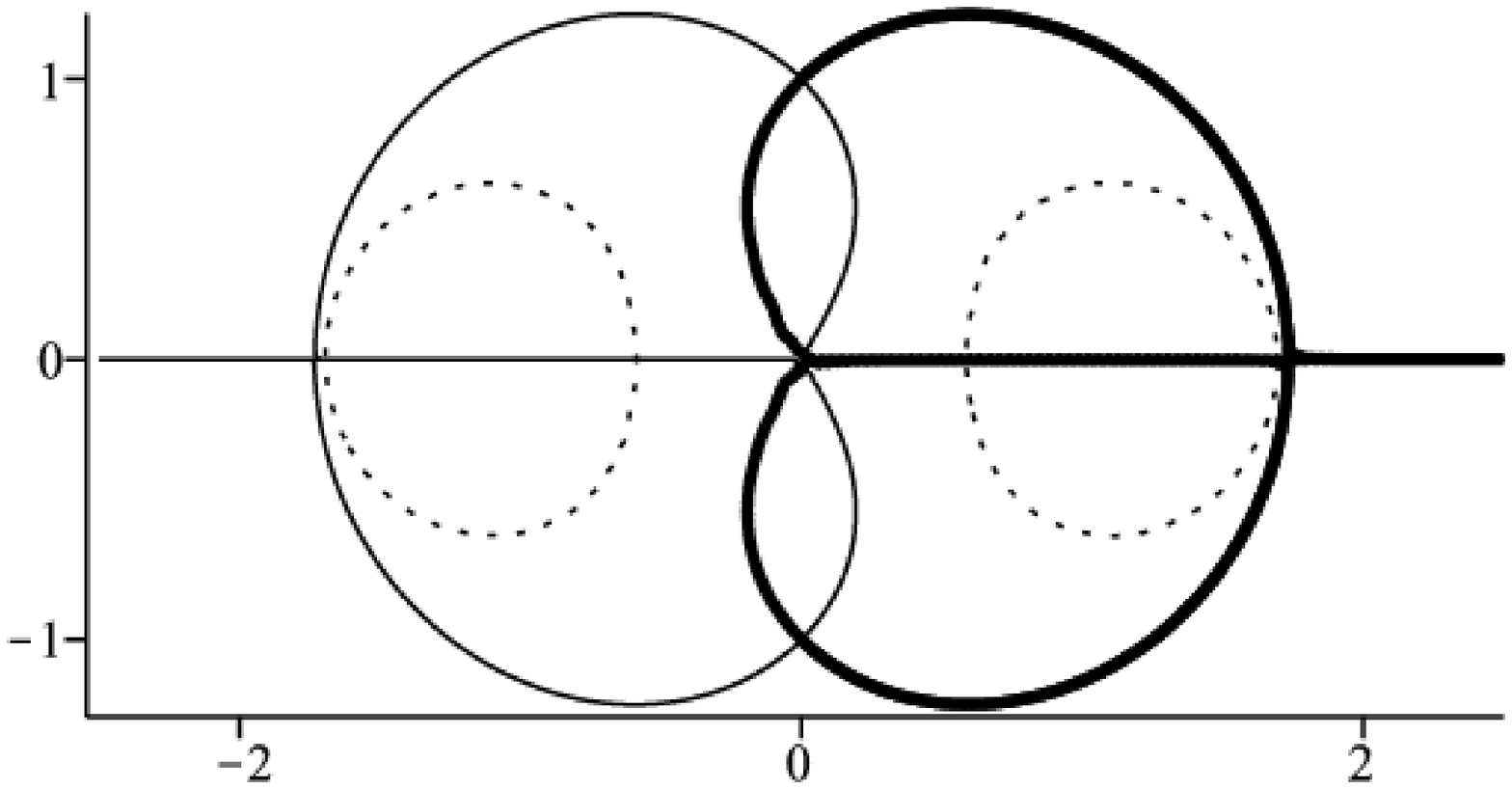}} 
\caption{Curves (\ref{I}) in dotted lines and (\ref{II}) 
in solid lines for $k=1.1$.} 
\label{f8} 
\end{figure}

\begin{figure} 
\epsfxsize=3.2in 
\centerline{\epsffile{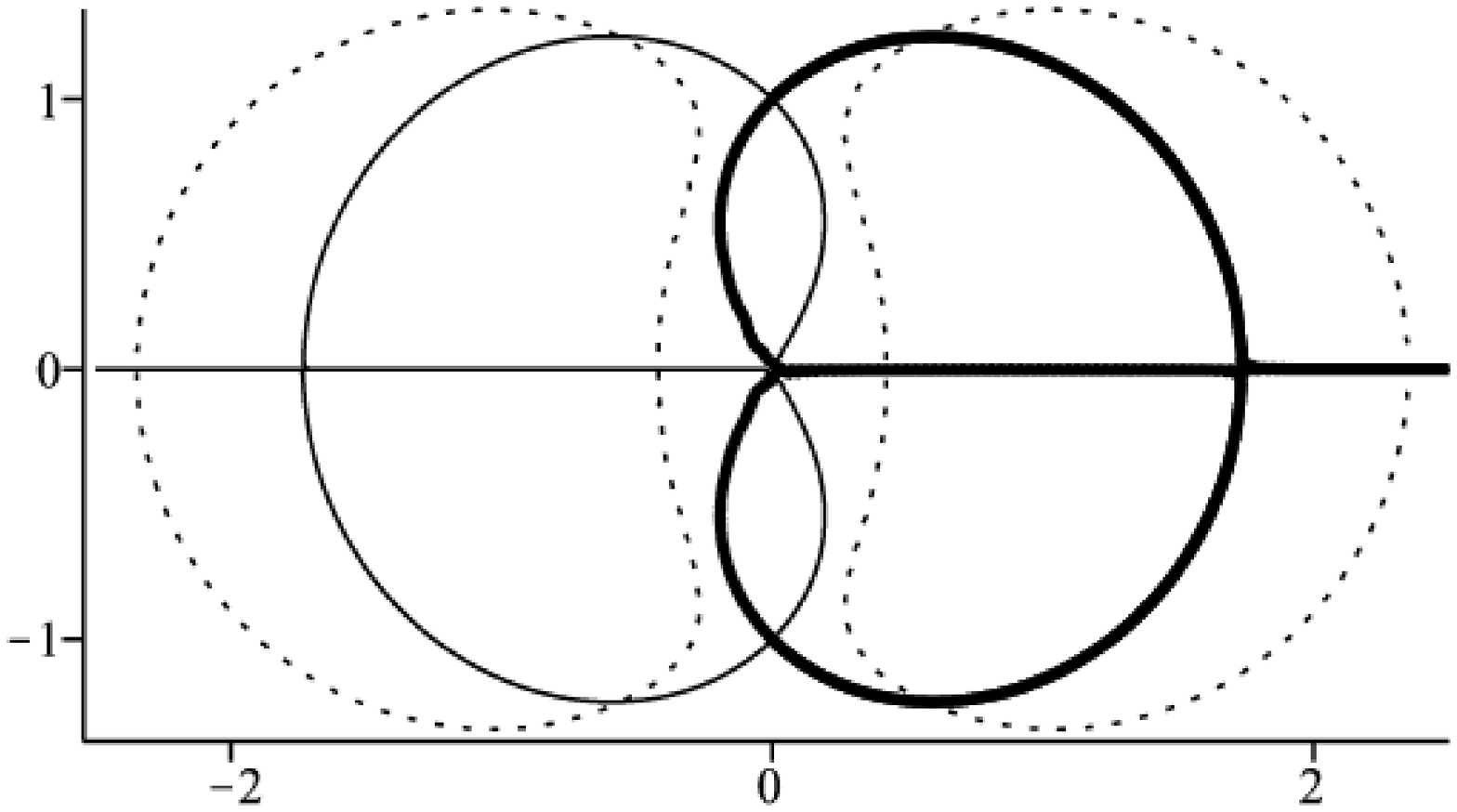}} 
\caption{Curves (\ref{I}) in dotted lines and (\ref{II}) 
in solid lines for $k=1.92$.} 
\label{f9} 
\end{figure}

\begin{figure} 
\epsfxsize=3.2in 
\centerline{\epsffile{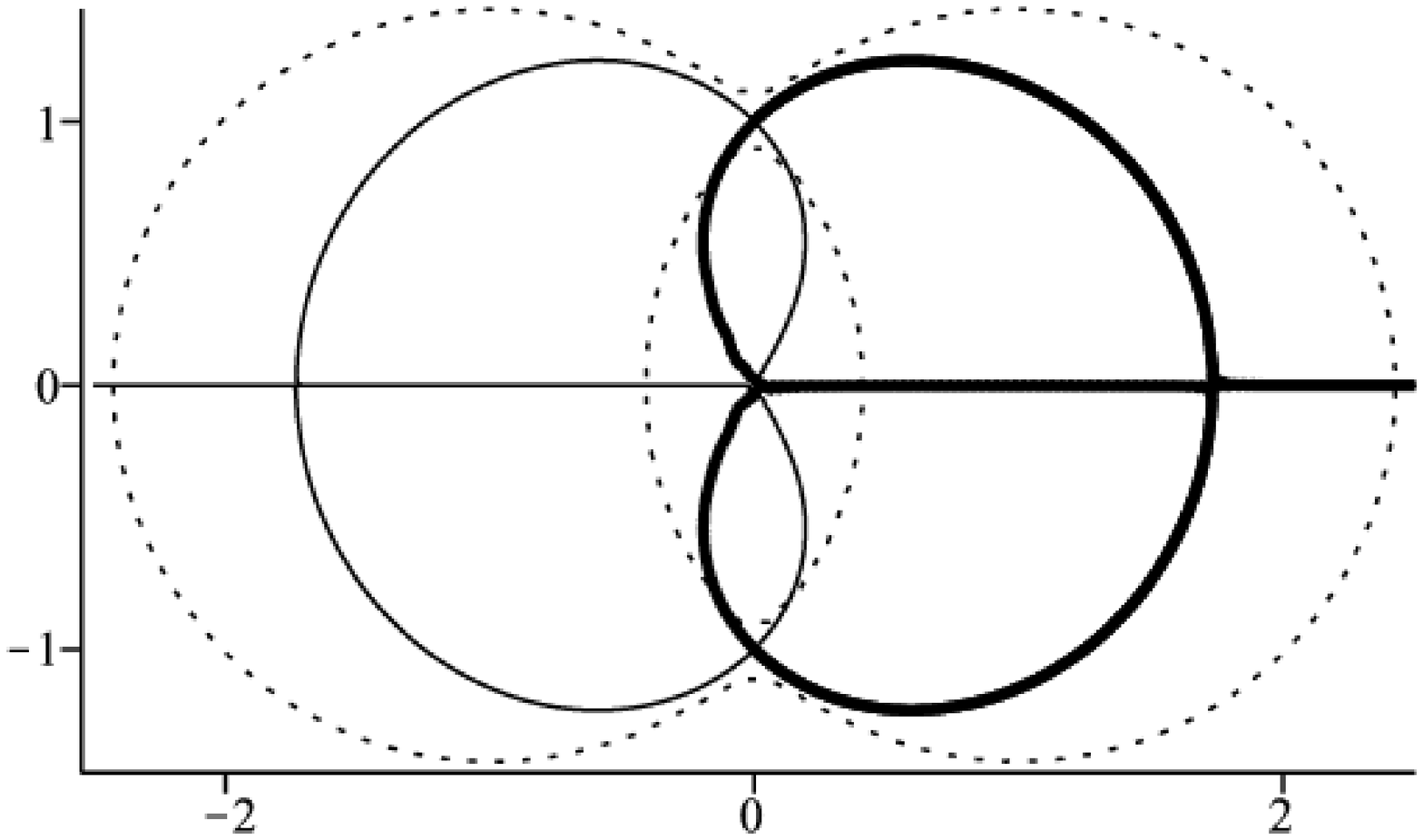}} 
\caption{Curves (\ref{I}) in dotted lines and (\ref{II}) 
in solid lines for $k=2.01$.} 
\label{f10} 
\end{figure}

Pictures of the algebraic curves 
defined by~(\ref{I}) and~(\ref{II})
are shown in 
Figures \ref{f8}--\ref{f10} for various values of the para\-meter~$k$. 
The part where 
\begin{equation}
\label{IIa}
\frac{(1+s^2)^2}{s^3}>0
\end{equation} 
is shown by a bold line. 
 
It is easy to see that there are always $4$ real solutions of~(\ref{I}): 
\begin{equation}\label{r} 
s=\pm \frac{k}{2}\pm\sqrt{\frac{k^2}{4}+1},
\end{equation} 
two of them in the interval $(-1,1)$ and two outside of 
this interval. 
 
After simplification and factoring out $s-\s$ 
from~(\ref{II}), we obtain 
\begin{equation}\label{15} 
s^2+\s^2=1-k^2|s|^2+|s|^4 
\end{equation} 
and 
\begin{equation}\label{16} 
s^2+\s^2=|s|^6-2|s|^4-|s|^2. 
\end{equation} 
Eliminating $s^2+\s^2$ from these two equations, we obtain 
\begin{equation}\label{14} 
p(r)=r^3-3r^2+(k^2-1)r-1=0,\quad\mbox{where}\quad r=|s|^2. 
\end{equation} 
The critical values of this polynomial $p$ 
are  
$$(k^2-4)\left(1\pm\frac{2}{9}\sqrt{12-3k^2}\right),$$ 
and they are both negative for $0<k<2$, both equal to $0$ for $k=2$
 and non-real for $k>2$.
As $p(0)<0$ we conclude that $p$ has exactly one positive root,
for all $k>0$.
We denote  this root by $r(k)$. 
 
The equation (\ref{15}) now gives 
\begin{equation}\label{s2s2bar}
s^2+\s^2=1-k^2r(k)+r^2(k),\quad\mbox{where}\quad r(k)=|s|^2,
\end{equation}
and this has solutions if and only if  
\begin{equation}\label{rkbar}
\left|1-k^2r(k)+r^2(k)\right|\leq 2 r(k).
\end{equation}
The equation $p(r(k))=0$ can be written in the form 
$$1-k^2 r(k)+r^2(k) =-r(k)-2r^2(k)+r^3(k)$$
and this yields
$$1-k^2 r(k)+r^2(k) 
=-2r(k)+r(k)\left(1-r(k)\right)^2\geq -2r(k).$$
Thus the absolute value sign can be dropped from~(\ref{rkbar})
and we obtain
$$1-k^2r(k)+r^2(k)\leq 2 r(k).$$ 
With $q(r)=r^2-(k^2+2)r+1$ we thus have  
$p(r(k))=0$ and $q(r(k))\leq 0$.  
Since $p(r)+q(r)=r^3-2r^2-3r=r(r+1)(r-3)$ 
we conclude that $r(k)\leq 3$. 
Hence $p(3)=3k^2-4\geq 0$. 

So
the equations~(\ref{15}) and~(\ref{16}) 
have common solutions if and only if 
$k\geq 2/\sqrt{3}$. 
If $s$ is a solution, then so are $-s$,  $\overline{s}$ and 
$-\overline{s}$. 
With $s=\sqrt{r(k)}e^{it}$ the equation~(\ref{s2s2bar})
takes the form 
\begin{equation}\label{cos2t}
\cos(2t)=\frac{1-k^2r(k)+r^2(k)}{2r(k)}.
\end{equation}
We find that 
if $k>2/\sqrt{3}$, $k\neq 2$, then
there are exactly $4$ solutions
of the system given by~(\ref{15}) and~(\ref{16}),
one in each open quadrant. 

We now determine for which solutions the inequality~(\ref{IIa}) is satisfied.
With $s=|s|e^{it}=\sqrt{r(k)}e^{it}$ we have
\begin{eqnarray*}
|s|^3\frac{(1+s^2)^2}{s^3}
&=&|s|^3\Rea \left(\frac{(1+s^2)^2}{s^3}\right)\\
&=&|s|^3\Rea \left(\frac{1}{s^3}+\frac{2}{s}+s\right)\\
&=&\cos(3t) +\left(2 r(k) +r^2(k)\right)\cos(t)\\
&=&\cos(t)\left(2\cos(2t)-1+2r(k) +r^2(k)\right).
\end{eqnarray*}
Using (\ref{cos2t}) and noting that $p(r(k))=0$
we obtain
\begin{eqnarray*}
|s|^3\frac{(1+s^2)^2}{s^3}
&=&\cos(t)\left(\frac{1-k^2r(k)+r^2(k)}{r(k)}-1+2r(k)+r^2(k)\right)\\
&=&\frac{\cos(t)}{r(k)}\left(1-(k^2+1)r(k)+3r(k)^2+r^3(k)\right)\\
&=&\frac{\cos(t)}{r(k)}\left(1-(k^2+1)r(k)+3r(k)^2+r^3(k)+p(r(k))\right)\\
&=&\frac{\cos(t)}{r(k)}\left(-2r(k)+2r^3(k)\right)\\
&=&2\cos(t)\left(r^2(k)-1\right).
\end{eqnarray*}
Thus (\ref{IIa}) is satisfied for the  solutions in the first 
and fourth quadrant
if $r(k)>1$ and for the solutions in the second and third quadrant 
if $r(k)<1$. Since $p(1)=k^2-4$ we see that the first alternative
holds for $k<2$ while the second one holds for $k>2$;
 cf.\ Figures~\ref{f9} and~\ref{f10}.

To summarize, 
we find that 
for $2/\sqrt{3}<k<2$
the inequality~(\ref{IIa}) is satisfied for the solutions in the 
first
and fourth quadrant.
Since the cosine is a proper map of degree $2$ from
$D^0$ onto the right half-plane, each of these solutions
corresponds to $2$ solutions for 
the original variable~$z$.
More precisely, the solution for $s$ in the first quadrant 
corresponds to two solutions for $z$ in the second and
fourth quadrant, and the solution for $s$ in the fourth
quadrant corresponds to two solutions for $z$ in the first and
third quadrant.

For $k>2$ the inequality~(\ref{IIa}) is satisfied for the solutions in the
second and third quadrant. Thus the preimages
of these solutions under the cosine are outside of the region~$D^0$.
For $k=2$ we find two solutions on the imaginary axis.
Again the preimages
under the cosine are outside of the region~$D^0$.

In addition to the solutions obtained from solving~(\ref{15}) 
and~(\ref{16}), we always have the $4$  
real solutions given by~(\ref{r}). 
Two of these solutions are positive, and these are the ones that
satisfy~(\ref{IIa}). 
Moreover, one of them is in the interval $(0,1)$ and one is in 
the interval $(1,\infty)$. In the original variable $z$ the 
first one corresponds to $2$ real solutions while
the second one corresponds to $2$ solutions on the imaginary 
axis.
 
Altogether we thus conclude  
that there are $4$ or $8$ points $z$ on $\gamma$ 
such that $f(z)$ is a cusp of $\Gamma$. 
Of the $4$ points corresponding to the 
real solutions~(\ref{r}), there is one on each coordinate 
semi-axis. 
We label these points as $z_1,z_2,z_3,z_4$ where 
$z_1>0,$ $z_2=ic,c>0$, $z_3=-z_1$ and $z_4=-z_2$. 
For $2/\sqrt{3}<k<2$ there are $4$ further solutions  
$w_1,w_2,w_3,w_4$ which we label such that 
$w_j$ is in the $j-$th 
quadrant; 
see Figure~\ref{f2}. 
The parameter $k=2$ corresponds to the limiting case
where $w_1,w_2,w_3,w_4\in\partial D^0$.
 
Now we determine the position of the cusps $f(z_j)$. 
Since $0<z_1<\pi/2$, we deduce from~(\ref{11})  
that $k\cos z_1/\sin^2 z_1=1$ and hence $k/\sin z_1=\tan z_1$. 
Thus  
$$f(z_1)=z_1-\frac{k}{\sin z_1}=z_1-\tan z_1<0.$$ 
Similarly we figure out where the other three cusps $f(z_j)$ 
are located and find that 
\begin{equation} 
\label{cusps} 
f(z_1)<0,\quad f(z_3)>0,\quad f(z_2)/i>0,\quad\mbox{and}\quad f(z_3)/i<0. 
\end{equation}

Thus we obtain the following result. 
\vspace{.1in} 
 
\noindent 
{\bf Lemma 1.} {\em The caustic $\Gamma$ has 
$4$ cusps if $0<k\leq 2/\sqrt{3}$ or $k\geq 2$ 
and it has $8$ cusps if $2/\sqrt{3}<k<2.$ 
For every $k>0$ there are $4$ cusps $f(z_j)$ on the 
coordinate axes located as in $(\ref{cusps})$. 
For $2/\sqrt{3}<k<2$ there are $4$ additional cusps 
$f(w_j)$.} 
\vspace {.1in} 
 
The tangent vectors at the cusps are horizontal 
on the real line and vertical on the imaginary line.  
This follows from the symmetry of $ 
\Gamma$ with respect to reflections in the coordinate axes. 
\vspace{.1in} 
 
\noindent 
{\bf Lemma 2.} {\em The tangent vector 
to $\Gamma=f(\gamma)$ is never  
vertical, except on the imaginary axis, 
and never horizontal, except on the real axis.} 
\vspace {.1in} 
 
{\em Proof.}  
By \cite[(2.4)]{Lyzzaik}
this tangent vector is  collinear to $\pm e^{it/2}$, 
where $t$ is the para\-meter defined in~(\ref{12}). 
So the tangent vector is horizontal or vertical if and only if
$g'(z)$ is real. 
Thus
$$k\frac{\cos z}{\sin^2 z}=k\frac{\cos z}{1-\cos^2z}=\pm 1.$$
This yields 
$$\cos z=\pm\frac{k}{2}\pm \sqrt{\frac{k^2}{4}+1}\in\R$$
which for $z\in D^0$ implies that $z$ is on 
the real or imaginary axis.
\vspace{.1in} 
 
A smooth curve will be called {\em convex} if its tangent 
vector turns to the left all the time. 
In other words, $\arg \zeta'(t)$ 
strictly increases, where $\zeta(t)$ is the para\-metrization 
of the curve.  
We will use the following fact: 
\vspace{.1in} 
 
\noindent 
{\bf Lemma 3.} \cite[Theorem 2.3]{Lyzzaik} {\em  
Each smooth piece of $\Gamma$, 
para\-metrized as explained above as $\zeta(t)=f(z(t))$ 
between the cusps, 
is convex. At the cusps the argument of the tangent  
vector jumps by $\pi$.} 
\vspace{.1in} 
 
Lemmas 1--3 are sufficient for the proof of all 
properties of $\Gamma$ we need. 
\vspace{.1in} 
\begin{center}  
{\bf 4. Proof of Propositions 1 and 2}  
\nopagebreak[4]  
\end{center}  
 
{\em Proof of Proposition 1.}
First we consider the case that $0<k\leq 2/\sqrt{3}$. The curve 
$\Gamma$ has $4$ cusps, one on each coordinate semi-axis. 
Begin tracing 
$\Gamma$ from the cusp $f(z_3)$ on the positive semi-axis, 
where its tangent has argument~$\pi$. 
As the argument of the tangent increases 
and can never reach $3\pi/2$, this arc ends 
at  the cusp on the negative ray of the imaginary axis. 
Both $\Rea\zeta$ and $\Ima\zeta$ are monotone on the 
arc, so it belongs to the $4$-th quadrant. 
The other three smooth arcs of $\Gamma$ are obtained 
by symmetry with respect to both axes. 
Thus $\Gamma$ is a simple Jordan curve, as shown in Figure~\ref{f4}, 
and the index of $\Gamma$ about any point $w$ 
can be only $0$ or~$\pm 1$.    
 
Now we consider the case that $2/\sqrt{3}<k<2$.  
Let $\gamma_0$ be the arc of $\gamma$ from $z_3$ to $z_4$, 
and put $\Gamma_0=f(\gamma_0)$. 
With $w_3\in\gamma_0$ defined as above,  
the points $f(z_3)$, $f(w_3)$ and $f(z_4)$ are consecutive 
cusps on $\Gamma_0$. Denote by $\Gamma_1$ the arc of $\Gamma_0$ 
from $f(z_3)$ to $f(w_3)$, and by $\Gamma_2$ the arc from 
$f(w_3)$ to $f(z_4)$. 
Then the tangent to~$\Gamma_1$ 
at $f(z_3)$ has argument $\pi$ and the argument increases but never 
reaches $3\pi/2$ on~$\Gamma_1$. 
It follows that $\Ima\zeta$ {\em decreases} on~$\Gamma_1$. 
 
Let $a\in (\pi,3\pi/2)$ 
be the argument of $\Gamma_1$ at $f(w_3)$. 
The next arc $\Gamma_2$ of $\Gamma_0$ begins 
at $f(w_3)$ with the argument of the tangent 
$a-\pi\in (0,\pi/2)$ and then the argument of 
the tangent increases but reaches the value $\pi/2$ 
only at the endpoint $f(z_4)$ on the negative 
imaginary axis. It follows that $\Ima\zeta$ {\em increases} 
on~$\Gamma_2$. 
So $\Gamma_0$ intersects every horizontal line at most twice. 
 
Another conclusion from these arguments is that $\Gamma_0$ 
belongs to the lower half-plane, except for the point $f(z_3)$. 
Indeed, $\Ima\zeta$ decreases on $\Gamma_1$ from $0$ 
to some negative value, and then $\Ima\zeta$ increases on 
$\Gamma_2$ and ends with a negative 
value. 
 
Let $\Gamma_3$ be the curve obtained by reflecting 
$\Gamma_0$ in the imaginary axis and 
reversing the orientation. Then the sum 
 $\Gamma_4=\Gamma_0+\Gamma_3$ intersects 
all horizontal lines at most $4$ times, 
and does not intersect 
horizontal lines in the upper half-plane. 
Let $\Gamma_5$ be the curve obtained by reflecting 
$\Gamma_4$ in the real axis and changing the orientation. 
Then the sum $\Gamma_6=\Gamma_4+\Gamma_5$ intersects 
every horizontal line at most $4$ times. 
On the other hand, $\Gamma_6=\Gamma$, and we conclude that 
the index of $\Gamma$ has absolute value at most~$2$.  
 
Finally we consider the case $k\geq 2$.
We shall actually assume that $k>2$ and make a remark about the 
changes for the case $k=2$ at the end.
Let $v_1$ and $v_2$ be the endpoints of $\gamma$ 
on the line $\Rea z=-\pi/2$ such that $\Ima v_2<\Ima v_1<0$; see  
Figure~\ref{f3}. 
Let $\gamma_0$ be the arc of $\partial D^-$ from $z_3$ to 
$z_4$. We break $\gamma_0$ into three arcs: 
\vspace{.1in} 
 
$\gamma_1$ from $z_3$ to $v_1$, 
 
$\gamma_2$ from $v_1$ to $v_2$, and 
 
$\gamma_3$ from $v_2$ to $z_4$. 
\vspace{.1in} 

\noindent
The images of $\gamma_j,\;0\leq j\leq 3$, under $f$ are 
denoted by $\Gamma_j$. 
 
The imaginary part $\Ima\zeta$ is decreasing on $\Gamma_1$ 
for the same reasons as before: 
$\Gamma_1$ is a convex curve which begins at $f(z_3)$ 
with horizontal tangent, and the argument of the tangent 
increases on $\Gamma_1$ but never reaches a vertical direction. 
 
The imaginary part $\Ima\zeta$ is decreasing on 
$\Gamma_2$ because  
$$\Ima f\left(-\frac{\pi}{2}-it\right)=-it.$$ 
 
Finally, $\Ima\zeta$ is monotone on $\Gamma_3$ 
because $\Gamma_3$ ends on the imaginary axis 
with a vertical tangent, and this tangent never 
else has a vertical direction. 
 
We conclude that $\Gamma_0$ intersects every horizontal line 
at most twice, 
 and $\Gamma_0$ belongs to the lower half-plane. 
As $\partial D^-$ is the union of $4$ curves obtained 
from $\Gamma_0$ by reflections in the axes, we conclude 
that $\partial D^-$ intersects each horizontal line 
at most $4$ times, 
so its index does not exceed $2$ by absolute value. 

The case $k=2$ corresponds to the case that $v_1=v_2$
so that the curve $\gamma_2$ degenerates to a point.
This case can be handled by the same argument.
This completes the proof of Proposition~1.

\vspace{.1in} 
  
{\em Proof of Proposition 2.} 
For $k>2$ the open set $D^+$ consists of $4$ regions 
$D_j^+$, $1\leq j\leq 4$, where $D_1^+$ intersects 
the positive real axis, $D_2^+$ intersects the 
negative real axis, $D_3^+$ intersects the positive 
imaginary 
axis and $D_4^+$ intersects the negative imaginary axis. 
 
Let $L=[\pi/2-it_0,\pi/2+it_0]=[v_1,-v_1]$ 
be the vertical segment of $\partial D_1^+$. 
It is easy to see that $t_0$ is defined by 
$$\sinh t_0=\frac{k}{2}-\sqrt{\frac{k^2}{4}-1}<1\quad\mbox{for}\quad k>2.$$  
An easy computation shows that $f(L)$ is the graph 
of a convex function $x=\phi(y)$. 
Indeed, we have 
$$\phi(y)=\frac{\pi}{2}-\frac{k}{\cosh y},$$ 
and this is convex for $\sinh y\in [-1,1]$. 
 
So $f(\partial D_1^+)$ consists of 
three convex curves. Two of them
are symmetric to each other
with respect to the real line and meet at a cusp
on the real line. As the variation
of the argument of the tangent to each of
these two curves is less than $\pi/2$,
the union of these curves is a graph of a function
$x=x_1(y)$.
The third curve is
a graph of a function $x=x_2(y)$.

It follows that every horizontal line intersects
$f(\partial D_1^+)$ at most twice,
so the index of $f(\partial D_1^+)$ with respect to
every point in the plane has absolute value at most~$1$.
The same argument works for $\partial D_2^+$. 
So $|I_w(f(\partial D_j^+))|\leq 1$ for every $w$ 
and $j=1,2.$ 
 
We now consider the curve $f(\partial D_4^+)$.
Let $\sigma_1$ be the arc between $z_4$ and~$v_2$.
Then $f(\sigma_1)$ is convex without a horizontal tangent
and thus intersects every horizontal line at most once.
Let $\sigma_2$ be the left boundary of~$D_4$; that is,
$\sigma_2=\{-\pi/2+it: t<-|v_2|\}$.
Since $f(-\pi/2+iy)=-\phi(y)+iy$ with the above function $\phi$
we see that  $f(\sigma_2)$
also intersects every horizontal line at most once.
Moreover, we find that $f(\sigma_1+\sigma_2)$ is smooth at 
the point $f(v_2)$. 
This implies that that $f(\sigma_1+\sigma_2)$ intersects every
horizontal line at most once.
By symmetry, the same holds for 
the image of the reflection of $\sigma_1+\sigma_2$
at the imaginary axis.
It follows that the index of 
the curve $f(\partial D_4^+)$ is at most~$1$.
By symmetry, the index of $f(\partial D_3^+)$ is also at most~$1$.

Finally the curves $f(\partial D_j^+)$, $j=3,4$ do 
not intersect because one of them belongs to the 
upper half-plane and the other one belongs to the lower half-plane. 
 
Thus $f(\partial D^+)$ is the union of $4$ 
curves, each of them has index of absolute value at most 
$1$ with respect to any point, 
and two of these curves belong to 
complementary half-planes. 
The conclusion follows for $k> 2$. 

If $k=2$, the domains $D_1^+$ and $D_2^+$ are not present.
The above arguments for $D_3^+$ and $D_4^+$ hold without change
and thus the conclusion of Proposition~2 follows also in this case.
\vspace{.2in} 
 
\begin{center}
{\bf 5. Proof of Proposition 3}
\end{center}
\vspace{.1in}

We begin with the following lemma.
\vspace{.1in}

\noindent
{\bf Lemma 4.} {\em Let $H_1,\ldots,H_n$ be closed half-planes
whose boundaries contain a point~$w$. Then there exist points
$w'\neq w$ arbitrarily close to $w$ which belong to at least
$[n/2]+1$ half-planes.}
\vspace{.1in}

{\em Proof.} 
Suppose that $w'\neq w$ belongs to the boundary of one of the 
half-planes $H_1,\ldots,H_n$. If $w'$ belongs to at most $[n/2]$
half-planes then $w''=2w-w'$ belongs to at least $[n/2]+1$
half-planes.
Since $|w''-w|=|w'-w|$, the conclusion follows.
\vspace{.1in}

{\em Proof of Proposition 3}.
Let $w$ be a point with maximal number $m$ of preimages.
We are going to prove that there is a neighborhood $W$
of $w$
such that every $w'\in W$ has also $m$ preimages.
Consider the full preimage $f^{-1}(w)=\{ z_1,\ldots,z_m\}$.
We say that $f$ is {\em locally surjective at $z_j$} if
the image
of every neighborhood of $z_j$ is a neighborhood of~$w$.
If $f$ is locally surjective at every $z_j$
then we are done.
By the implicit function theorem $f$ is locally surjective at every $z_j$
which is not on the critical curve.

For points on the critical curves we shall use the classification
given by Lyzzaik~\cite[Definition~2.2]{Lyzzaik}
into points of the first, second and third kind.
We refer to his paper
for the precise definition, but note that
for the function $f$ defined by~(\ref{deff})
which was considered in Propositions 1 and 2 the points
of the first kind are the points $z$ on the critical curve for
which $f(z)$ is a cusp of the caustic while all other points on the 
critical curve are of the second kind. Thus points of the third kind
do not occur for that function, but they may occur in the more general
situation considered now.
It follows from Lyzzaik's local description
of harmonic maps~\cite[Theorem 5.1~(b)]{Lyzzaik} 
that $f$ is 
locally surjective at points of the first kind.

Suppose now that there is a point $z=z_j\in f^{-1}(w)$ where
$f$ is not locally surjective. Then $z$ is of the second or
third kind and Lyzzaik's 
results~\cite[Theorem 5.1~(a), Proposition~6.1]{Lyzzaik}
imply that there are topological
discs $U$ containing $z$ and $V$ containing $w$ with
the following properties.

Let $\gamma$ be the part of the critical curve in $U$ and
$\Gamma=f(\gamma)$. Then $\Gamma$ is a union of simple
curves
$\Gamma_k$ from $w$ to $\partial V$ which are smooth,
convex and have a common tangent at~$w$.
These curves $\Gamma_k$ split $V$ into some open curvilinear
sectors $D_k$. Each of these sectors is either completely
covered by the image
$f\vert_U$ or is disjoint from this image.
Let $g=f\vert_U$. We assume without loss of generality
that $f(\zeta)\neq w$ for $\zeta\in
\overline U\backslash\{ z\}$.
\vspace{.1in}

\noindent
{\bf Lemma 5.} {\em There can be at most one sector $D_k$
which is disjoint
from $g(U)$. The closure of this sector is
contained in a set of the form $(H^0\cap V)\cup\{ w\}$,
where $H^0$ is an open half-plane whose boundary contains
$w$.
If such a sector $D_k$ indeed exists then
all points in the other sectors have at least two
preimages under~$g$.}
\vspace{.1in}

{\em Proof.} To prove the first statement, we suppose to the
contrary that there are two sectors $D_k$ and $D_l$
disjoint from $g(U)$. We consider a closed disc $B$
around $w$ which is disjoint from $f(\partial U)$.
Such a disc exists because $f(\partial U)$ is a
compact set that
does not contain~$w$.
The preimage $K$ of this disc $B$ in $U$
is compact and
connected. (It is connected because every component
has to contain
a preimage of~$w$, and we assume that there is only one
such preimage in $\overline{U}$).
Moreover, $K$ is a neighborhood of $z$
because $g$ is continuous.
Now the set $B\backslash(D_k\cup D_l\cup\{ w\})$
is disconnected while its preimage coincides with
the preimage of $B\backslash\{ w\}$ that is equal to
$K\backslash\{ z\}$, and this set is connected because $K$
is a neighborhood of~$z$.
This contradiction proves the first statement.

To prove the second statement, we assume that $V$
is a round disc,
which does not restrict generality. 

We say that a region $D$ is {\em locally convex}
at a point
$\zeta\in\partial D$ if the intersection of $D$ with a
sufficiently small round disc centered at $\zeta$ is strictly
convex. Local convexity at each boundary point implied
strict convexity of the region.

We claim that a sector
$D_k$ which is disjoint from $g(U)$ must be strictly convex.
Indeed, it is bounded by two strictly
convex curves $\Gamma_k,\Gamma_{k+1}$ with a
common tangent at $w$
and an arc of the circle $\partial V$.
To see that the curves $\Gamma_k,\Gamma_{k+1}$
are convex in the ``right direction'', consider
a para\-metrization $\Gamma_k(t)$. It follows from the local
description
of the ``folds'' in \cite{Lyzzaik} that the
number of preimages
of a point {\em decreases by $2$} as we cross $\Gamma_k$
in the direction of the normal $\Gamma_k^{\prime\prime}$,
so $\Gamma_k^{\prime\prime}$ points to the inside of $D_k$.

So $D_k$ is locally convex at every point, except possibly
at~$w$. To prove the local convexity at~$w$, we have to exclude 
the possibility that $w$ is a cusp of $\partial D_k$.
But this possibility is excluded by the results of Lyzzaik quoted
before Lemma 5 which says that $f$ is locally surjective 
at points of the first kind.

However, there is also a short independent argument showing
that $w$ cannot be a cusp:
if $\partial D_k$ is locally convex
at every point except $w$ and is not locally convex at $w$
and thus has a cusp at~$w$, then there exists an open half-plane
$H$ such that $D_k\cup H\cup\{ w\}$ is a full neighborhood
of $w$. Let $\phi(\zeta)=a\zeta+b$ be an affine function
such that $\phi(w)=0$ and $\Rea \phi(\zeta)>0,\;\zeta\in H$.
Then the composition $\Rea \phi\circ g$ is a non-constant
harmonic
function having a minimum at~$z$, which is impossible.

This proves that $D_k$ is locally convex at~$w$,
and thus $D_k$ is strictly convex.
So we can take the tangent line to $D_k$ at $w$ as the
boundary of $H^0$ and the second statement of the lemma will
hold.

To prove the third statement we notice that the number
of preimages
changes by an even number when the point $w'$
crosses the caustic.

This completes the proof of Lemma~5.
\vspace{.1in}

Now we can complete the proof of Proposition~5.
Let $z_1,\ldots,z_n$ be the points in $f^{-1}(w)$
where $f$ is not locally surjective,
and $z_{n+1},\ldots,z_m$ be the points
in $f^{-1}(w)$ where
$f$ is locally surjective. According to Lemma~5,
to each $z_j$ with $1\leq j\leq n$ we can
associate a closed half-plane $H_j$
having $w$ on the boundary so that each point
$w'\in H_j\backslash\{ w\}$
which is close enough to $w$ has at least two preimages
 close to~$z_j$. By Lemma 4 there exists $w'\neq w$
which belongs to at least $[n/2]+1$ of these
these half-planes.
This point has at least $2([n/2]+1)+m-n$ preimages.
Since $2([n/2]+1)>n$ and thus $2([n/2]+1)+m-n>m$ for $n\geq 1$
and since $w$ was assumed to have the maximal number $m$
of preimages, we conclude that $n=0$.
Thus $f$ is
locally surjective at all preimages of the point~$w$.
This completes the proof.
\vspace{.1in}

We finish this section with a short outline of an example
communicated to us by A. Gabrielov
showing that Proposition~3 does not hold for general smooth
maps. Begin with a smooth map of the unit disc onto
the sector $\{ z=x+iy:|z|<1,x\geq 0,y\geq 0\}$ given by
$(x,y)\mapsto (x^2,y^2)$. Composing
this with a smooth homeomorphism we obtain a smooth map of
the unit disc onto the sector
$\{ z=x+iy:|z|<1,|\arg z|\leq\pi/(2m)\}$, where $m$
is an integer.
Let $B_1,\ldots,B_m$ be discs with pairwise disjoint closures 
in the plane. We define a smooth map $f$ in these discs
so that the discs are mapped onto $m$ disjoint sectors in 
the unit disc with common vertex at $0$ and of opening
$\pi/m$. Then we extend our map to the complement of
the discs $B_j$ so that the resulting map is smooth.
It is easy to see that this construction can be performed
so that for the resulting map each point except $0$
has at most $5$ preimages. Thus we obtain a smooth map
for which every point except $0$ has at most $5$ preimages
while $0$ has at least $m$ preimages, where $m$ is
arbitrarily large.
\vspace{.2in}

\begin{center} 
{\bf 6. Remarks} 
\end{center} 
\vspace{.1in}
 
1. It follows from our proof that the number of solutions of
equation~(\ref{1}) is constant in the complementary components of
the caustic and, given $k$ and a value $w$ which is not on the caustic,
the number of solutions can be computed from the index of the
caustic with respect to~$w$.
Since the caustic is an explicitly given curve,
the number of solutions can actually be computed.

To demonstrate that any number of solutions
between~$1$ and~$6$ can actually occur, 
we pick appropriate values of $w$ 
from our Figures \ref{f4}--\ref{f7a}, or similar figures 
for other values of~$k$. 
 
For example, $k=1.92$ and $w=0.67i$ gives $6$ solutions 
$1.5363458i,$ 
$-0.9885626i,$ 
$\pm 1.2603941+0.9732810i,$ 
$\pm 1.4617539+0.7738876i,$ 
 
\vspace{.1in} 
 
2. As mentioned at the end of section 2,
the region with $6$ solutions exists 
for $k\in (2/\sqrt{3}, k_0),$ 
where 
$$k_0=\frac{\pi^2}{2\sqrt{\pi^2-4}}\approx 2.0368$$
is determined from the equations 
$$\left|\frac{\cos z}{\sin^2 z}\right|=\frac{1}{k},\quad \Rea(z-k/\sin z)=0,\quad
z=\frac{\pi}{2}-it.$$
This is the condition that $\Rea f(v_1)=0,$ where $v_1$ 
was defined in the proof of Proposition 1 for $k>2$ in 
Section 4. 
\vspace{.1in} 
 
3. Equation (\ref{shear}) can be rewritten in a form similar to 
equation~(\ref{1}).  
Suppose that $\alpha\notin\{1,-1\}$ and put 
$u=z-\alpha\overline{z}.$ 
Then $\overline{z}=(\overline{\alpha}u+\overline{u})/(1-|\alpha|^2),$ 
and the equation becomes 
$$u=\arcsin\frac{k_1}{\overline{\alpha}u+\overline{u}+\overline{w_1}},$$ 
where $k_1=k(1-|\alpha|^2)$, and $w_1=w(1-|\alpha|^2)$. 
Now we take the sine on both sides and conjugate to obtain 
\begin{equation}\label{shear2} 
u=\frac{k_1}{\sin\overline{u}}-\alpha\overline{u}-w_1. 
\end{equation} 
This can be considered as a perturbation of the equation~(\ref{1}) 
by the term~$\alpha\overline{u}$. 
Orientation-preserving solutions of equation~(\ref{shear2}) 
are attracting fixed points of the anti-analytic entire function 
$h(u)=k_1/\sin\overline{u}-\alpha\overline{u}-w_1$. 
Fatou's theorem says that every attracting  fixed point attracts a trajectory 
of a singular value. If $\alpha=0$ the function has $3$ singular values, 
so there are at most $3$ attracting fixed points. 
This is the crucial part of the argument in~\cite{HL}. 
For $\alpha\neq 0$, 
the function $h$ has infinitely many critical values, 
so the dynamical proof breaks down at 
this point. 
 
\begin{figure}[t] 
\epsfxsize=4.5in 
\centerline{\epsffile{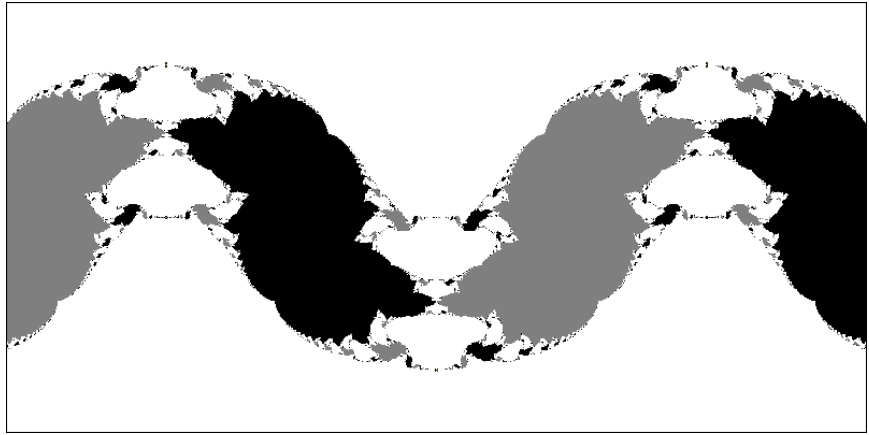}} 
\caption{Basins of attraction of  
$h(z)=1.92/\sin\overline{z}+0.67i$.  
} 
\label{f11} 
\end{figure}

In the case $\alpha=0$ considered in this paper, $h$ indeed can have $3$ 
distinct attracting points. This happens for those values 
of the para\-meters $k$ and $w$  
for which we have $3$ orientation-preserving solutions, 
for instance for $(k,w)=(1.92, 0.67i)$ as in the example above. 
Figure~\ref{f11} shows the partition of the plane into 
three domains of attraction of the fixed points: 
the attracting basin of $1.5363458i$ is shown in white, 
that of $1.4617539+0.7738876i$ in black and that of 
$-1.4617539+0.7738876i$ in gray. 
We mention that the Fatou set of this function is the union
of these attracting basins. Since this function has no wandering
domains~\cite{BKL} and no Baker domains~\cite[Theorem~A]{RS},
this can be deduced from relations between singularities
of the inverse and periodic components~\cite[Theorem~7]{Bergw}.
The Julia set of this function has zero area~\cite[Theorem~3]{Jank}, 
and it is not visible in the pictures.

{\em W. B.: Mathematisches Seminar, 
Christian-Albrechts-Uni\-ver\-si\-t\"at zu Kiel, 
Lude\-wig-Meyn-Str.~4, 
D-24098 Kiel, 
Germany 
 
bergweiler@math.uni-kiel.de 
 
A. E.: Department of Mathematics, 
Purdue University, West Lafayette, IN 47907, USA 
 
eremenko@math.purdue.edu} 
 
\newpage 
 

\end{document}